 \documentclass[twoside,reqno,11pt]{fcaa}

\usepackage{graphicx}
\usepackage{epsfig}
\usepackage{amsthm}
\usepackage{amsmath}
\usepackage{latexsym}
\usepackage{amsfonts}
\usepackage{amssymb}

\newtheoremstyle{theorem}
  {15pt}          
  {15pt}  
  {\sl}  
  {\parindent}
  {\sc}  
  {. }   
  { }    
  {}     
\theoremstyle{theorem}

\newtheoremstyle{defi}
  {15pt}          
  {15pt}  
  {\rm}  
  {\parindent}     
  {\sc}  
  {. }    
  { }    
  {}     
\theoremstyle{defi}

 \def\proofname{\indent\rm P\ r\ o\ o\ f.}
 
 \def\qed{\hfill$\Box$} 

 \usepackage{hyperref} 

 \def\theequation{\arabic{section}.\arabic{equation}}

 \def\themycopyrightfootnote{\vspace*{3pt}
 \copyright \, Year\,  Diogenes  Co., Sofia
   \par  \noindent pp. xxx--xxx, DOI: ..................
   \hfill  \vspace*{-36pt}
  }

 \title[Generalized fractional Poisson distribution]{Generalization of the fractional Poisson distribution \\ [4pt]}
 \author[R. Herrmann]{Richard Herrmann $^1$}

\begin{document}

 \vbox to 2.5cm { \vfill }


 \bigskip \medskip

 \begin{abstract}

A generalization of the Poisson distribution based on the generalized Mittag-Leffler function $E_{\alpha, \beta}(\lambda)$ is proposed and the raw moments are calculated algebraically in terms of Bell polynomials. It is demonstrated, that the proposed distribution function contains the standard fractional Poisson distribution as a subset. A possible interpretation of  the additional parameter $\beta$ is suggested.
 \medskip


 \smallskip


 \end{abstract}

 \maketitle

 \vspace*{-16pt}



 \section{Introduction}\label{sec:1}
Since its early beginnings, fractional calculus has developed as an ever growing powerful tool to model complex
phenomena in different branches of  e.g. physics or biology. 

A major field of research is the area of fractional Brownian motion and anomalous diffusion phenomena. Closely related are random walk problems and the investigation of stochastic processes,
where e.g. survival probabilities $R(t)$ differ from the exponential law, which is known from the classical Poisson process.  

Based on the concept of non-exponential power law Markov-processes, the fractional extension of the classical Poisson process leads to the standard fractional Poisson distribution function of the form
\cite{mee99},\cite{pod99},\cite{met00},\cite{rep00},\cite{las03}
\begin{eqnarray}
\label{poistart}
p_k (\alpha,\lambda)&=&\frac{(\nu \lambda^{\alpha })^{k}}{k!}\sum\limits_{n=0}^{\infty }\frac{(k+n)!}{n!}\frac{(-\nu \lambda^{\alpha })^{n}}{\Gamma (\alpha (k+n)+1)} \\
&=&   \frac{(\nu \lambda^{\alpha })^{k}}{k!} \frac{d^k}{d\lambda^k} E_\alpha(\lambda)|_{-\nu \lambda^{\alpha }}
\\
&=&   (-1)^k \big( \frac{\lambda^{k}}{k!} \frac{d^k}{d\lambda^k} E_\alpha(\lambda)\big) |_{-\nu \lambda^{\alpha }}
\qquad 0<\alpha\leq 1
\end{eqnarray}
which may be considered as a generalization of the classical Poisson distribution represented by the introduction of the Mittag-Leffler function $E_\alpha(\lambda)$ and its derivatives \cite{ml}.

The mean $\mu$ and the variance $\sigma$ for this fractional Poisson distribution follow  as \cite{las03}
\begin{eqnarray}
\label{poimoments1}
\mu&=& \frac{\nu \lambda^{\alpha }}{\Gamma (1+\alpha)} \\
\label{poimoments2}
\sigma &=&  \mu + \mu^2
( \frac{\sqrt{\pi}\,\Gamma(1+a)}{2^{2 \alpha - 1}\Gamma(1/2 + \alpha)} -1)
\end{eqnarray}

The distribution function (\ref{poistart}) is the result of the underlying process and allows for a direct comparison of experimental data with the corresponding model assumptions.

In the following, we will consider the inverse problem:
We will present fractional generalizations of distribution functions, based on generalized Mittag-Leffler functions. 

The motivation for this approach is straight forward: A given experimental setup will generate distribution data first, which may be directly compared with generalized distribution functions. The underlying microscopic process may be 
investigated in a second step.

A large amount of literature covers extensions of the standard Mittag-Leffler function by introducing additional parameters and
functional behavior (see e.g. \cite{kir99}, \cite{kil13} and references therein), but a  direct e.g. physical  interpretation of the additional parameters is a still open question in many cases. 

In the following we therefore will apply an approach, which allows for a wide class of Mittag-Leffler type functions to determine the corresponding distribution functions.  This distributions may be directly compared with observed phenomena, which will turn out to be helpful for a broader understanding of the functionality and scope of possible applications.

\section{Nomenclature}\label{sec:2}

\setcounter{section}{2}
\setcounter{equation}{0}
\setcounter{theorem}{0}

In order to calculate the major properties of a generalized Poisson distribution based on the generalized Mittag-Leffler function, let us first recall the major properties of the classical Poisson distribution, which is given by:  
\begin{equation}
\label{poipn}
p_k(\lambda) = \frac{1}{N}\frac{\lambda^k}{k!}
\end{equation}
with the normalization constant $N$, which is determined by the requirement of normalizability of the distribution, which coincides with the zeroth raw moment $\mu_0$: 
\begin{eqnarray}
\label{poinorm}
\mu_0 = \sum_{k=0}^{\infty}p_k(\lambda) = \frac{1}{N}\sum_{k=0}^{\infty}\frac{\lambda^k}{k!}=  \frac{1}{N} e^{\lambda} = 1
\end{eqnarray}
Higher raw moments are then given as:
\begin{eqnarray}
\label{poimoments}
\mu_n &=&  \sum_{k=0}^{\infty} k^n p_k(\lambda) = e^{-\lambda}\sum_{k=0}^{\infty}k^n \frac{\lambda^k}{k!} \equiv
B_n(\lambda) 
\end{eqnarray}
where (\ref{poimoments}) is the defining equation for the  Bell polynomials $B_n$, which obey  the recursion relation\cite{bel34}
\begin{eqnarray}
\label{poinbellrec}
\frac{d}{dx}  B_n(x) &=& \frac{1}{x}B_{n+1}(x) - B_n(x) 
\end{eqnarray}
and are related to 
the Stirling numbers of second kind $S(n,k)$ \cite{sti30} via Dobrinski's formula\cite{dob77}:
\begin{eqnarray}
\label{poister}
B_n(\lambda) &=& \sum_{k=1}^{n} S(n,k) \lambda^k
\end{eqnarray}

The first raw moments result as:
\begin{eqnarray}
\label{poimoments12}
\mu_1 &=&  \lambda\\
\mu_2 &=& \lambda (\lambda+1)
\end{eqnarray}
and therefore the mean $\mu$ and variance $\sigma$ follow as centered moments:
\begin{eqnarray}
\label{poimoments12}
\mu &\equiv& \mu_1 =  \lambda\\
\sigma &\equiv& \mu_2 - \mu_1^2= \lambda (\lambda+1) -\lambda^2 = \lambda
\end{eqnarray}

\section{Generalized fractional Poisson distribution}\label{sec:3}

\setcounter{section}{3}
\setcounter{equation}{0}
\setcounter{theorem}{0}
We now extend the definition of the classical Poisson distribution $p_k(\lambda)$ by applying the canonical fractionalization procedure: 

We replace the faculty in (\ref{poipn}) by the $\Gamma$ function, introduce as an example  two additional parameters $\alpha, \beta $  and define a generalized fractional Poisson distribution $p_k(\lambda, \alpha,\beta)$ via:
\begin{equation}
\label{poipnM}
p_k(\lambda, \alpha,\beta) = \frac{1}{N}\frac{\lambda^k}{\Gamma(\beta + \alpha k)}
\end{equation}
the normalization constant $N$ follows from the requirement of normalizability of the distribution:
\begin{eqnarray}
\label{poinormM}
\mu_0 = \sum_{k=0}^{\infty}p_k(\lambda, \alpha,\beta) = \frac{1}{N}\sum_{k=0}^{\infty}\frac{\lambda^k}{\Gamma(\beta + \alpha k)}=  \frac{1}{N} E_{\alpha, \beta}(\lambda) = 1
\end{eqnarray}
with the generalized Mittag-Leffler function  $E_{\alpha, \beta}(\lambda)$  which is defined as\cite{wi}:
\begin{eqnarray}
\label{poiML}
E_{\alpha, \beta}(\lambda) = \sum_{k=0}^{\infty}\frac{\lambda^k}{\Gamma(\beta + \alpha k)}
\end{eqnarray}
The distribution is therefore given by:
\begin{equation}
\label{poipnM}
p_k(\lambda, \alpha,\beta) = \frac{1}{E_{\alpha, \beta}(\lambda)}\frac{\lambda^k}{\Gamma(\beta + \alpha k)}
\end{equation}
Higher raw moments may be calculated by an iterative procedure, which we sketch for the case $n=1$. 

Starting with unnormalized zeroth moment, we apply the operator $\frac{d}{d \lambda} \lambda$:
\begin{eqnarray}
\label{poimoment1}
E_{\alpha, \beta}(\lambda) &=&  
\sum_{k=0}^{\infty} \frac{\lambda^k}{\Gamma(\beta + \alpha k)} \\
\frac{d}{d \lambda} \lambda\cdot E_{\alpha, \beta}(\lambda) &=&
  \sum_{k=0}^{\infty} (k+1) \frac{\lambda^{k}}{\Gamma(\beta + \alpha k)} 
\end{eqnarray}
Using the commutation relation
\begin{eqnarray}
\label{poicomm}
[ \frac{d}{d\lambda},\lambda] &=& 1 
\end{eqnarray}
the result for the first raw moment follows as:
 \begin{eqnarray}
\lambda \frac{d}{d \lambda} \cdot E_{\alpha, \beta}(\lambda)+ E_{\alpha, \beta}(\lambda)   &=& 
\sum_{k=0}^{\infty} k \frac{\lambda^{k}}{\Gamma(\beta + \alpha k)} +  \sum_{k=0}^{\infty}  \frac{\lambda^{k}}{\Gamma(\beta + \alpha k)} \\
\lambda \frac{d}{d \lambda}\cdot E_{\alpha, \beta}(\lambda)  &=& 
 \sum_{k=0}^{\infty} k \frac{\lambda^{k}}{\Gamma(\beta + \alpha k)} 
\end{eqnarray}
This procedure may be applied  iteratively $n$ times to obtain the n-th raw moment.

We obtain an iterative solution for the raw moments:
 \begin{eqnarray}
\mu_n  &=& \frac{1}{E_{\alpha, \beta}(\lambda)}
(\lambda \frac{d}{d \lambda})^n \cdot E_{\alpha, \beta}(\lambda)  
\end{eqnarray}
which results in a n-fold application of the Euler-operator $J_E = x \frac{d}{dx}$.

In order to deduce a closed form solution, we perform a normal ordering of the n-fold Euler operator  $J_E^n$:

With the settings
\begin{eqnarray}
\label{poimoment1t}
x^n &\rightarrow&  \lambda^n E_{\alpha, \beta}^{(n)}(\lambda)\\
 &=&  (:\lambda \frac{d}{d \lambda}:)^n E_{\alpha, \beta}(\lambda)
\end{eqnarray}
where $E_{\alpha, \beta}^{(n)}(\lambda)$ denotes the n-th derivative of the Mittag-Leffler function with respect to $\lambda$
and (::) indicates the normal ordered product \cite{wic50}
\begin{eqnarray}
\label{poinormal}
(:\lambda \frac{d}{d \lambda}:)^n = \lambda^n \frac{d^n}{d \lambda^n}
\end{eqnarray}
the iteration scheme is now  isomorphic to the recurrence relation (\ref{poinbellrec}) and we  obtain for the raw moments
\begin{eqnarray}
\label{poiBell}
\mu_n &=& \frac{1}{E_{\alpha, \beta}(\lambda)} B_n(:\lambda \frac{d}{d \lambda}:)\cdot  E_{\alpha, \beta}(\lambda)\\
          &=& \frac{1}{E_{\alpha, \beta}(\lambda)} \sum_{k=1}^{k=n} S(n,k)(\lambda^n E_{\alpha, \beta}^{(n)}(\lambda))
\end{eqnarray}
\begin{figure}
\begin{center}
\includegraphics[width=11.6cm]{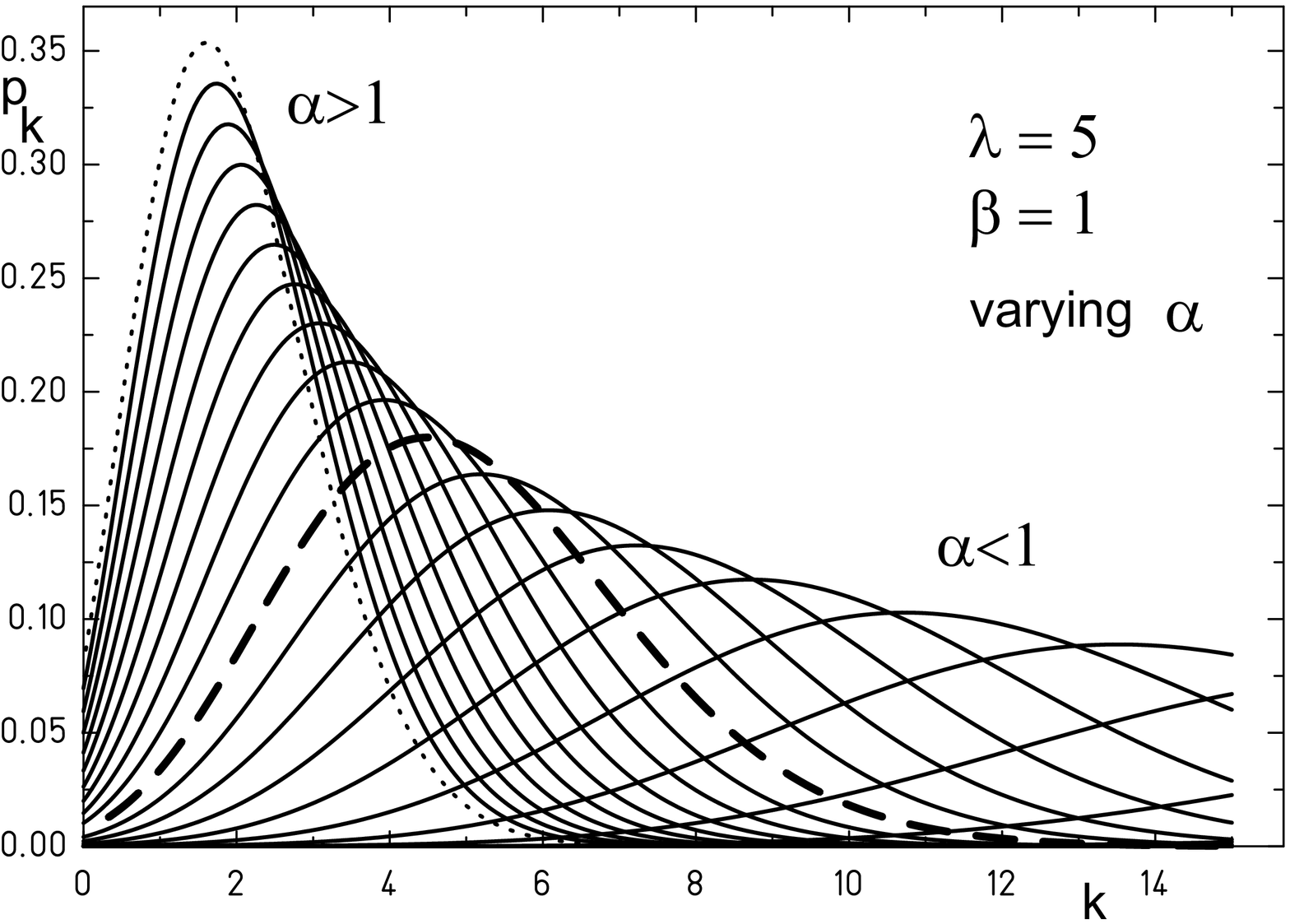}    
\caption{For $\lambda=5$ and $\beta = 1$ the probability distribution $p_k(\lambda=5, \alpha, \beta=1)$ from  (\ref{poipnM}) is plotted for $1.5 \geq \alpha \geq 0.5$ in $0.05$ steps. Thick dashed line indicates $\alpha=1$, which corresponds to the classical Poisson distribution. 
}
\label{fig1}
\end{center}
\end{figure}
\begin{figure}
\begin{center}
\includegraphics[width=11.6cm]{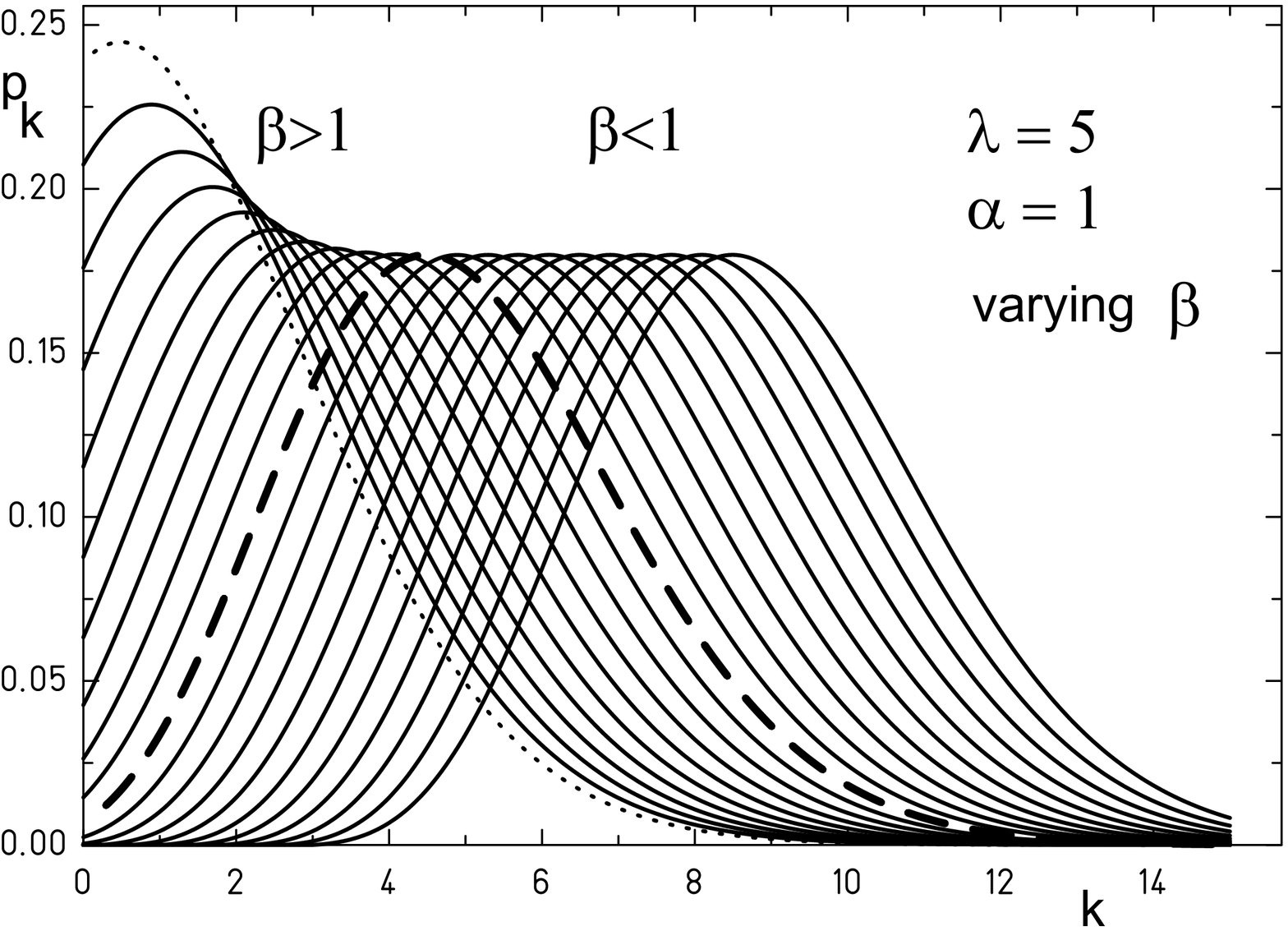}    
\caption{For $\lambda=5$ and $\alpha = 1$ the probability distribution $p_k(\lambda=5, \alpha=1, \beta)$ from  (\ref{poipnM}) is plotted for $4 \geq \beta \geq -4$ in $0.4$ steps. Thick dashed line indicates $\beta=1$, which corresponds to the classical Poisson distribution. 
}
\label{fig2}
\end{center}
\end{figure}
%
The mean $\mu$ and variance $\sigma$ follow as the lowest centered moments:
\begin{eqnarray}
\label{poimoments12x}
\mu \equiv \mu_1  &=&  \frac{1}{E_{\alpha, \beta}(\lambda)} \lambda  E^{(1)}_{\alpha, \beta}(\lambda)\\
&=&  \lambda \frac{d}{d \lambda} \log (E_{\alpha, \beta}(\lambda))\\
\label{poimoments12y}
\sigma \equiv \mu_2 - \mu_1^2 &=&
\frac{1}{E_{\alpha, \beta}(\lambda)}(
 \lambda^2  E^{(2)}_{\alpha, \beta}(\lambda)+
 \lambda  E^{(1)}_{\alpha, \beta}(\lambda))
 -\mu_1^2\\
 &=&  ( \lambda \frac{d}{d \lambda} \log (E^{(1)}_{\alpha, \beta}(\lambda)))
( \lambda \frac{d}{d \lambda} \log (E_{\alpha, \beta}(\lambda))) \nonumber \\
&&-( \lambda \frac{d}{d \lambda} \log (E_{\alpha, \beta}(\lambda)))^2 
+\lambda \frac{d}{d \lambda} \log (E_{\alpha, \beta}(\lambda))
 \end{eqnarray}
Optionally, with 
\begin{eqnarray}
\label{poiderivE}
\frac{d}{dx}E_{\alpha, \beta}(x) = \frac{1}{\alpha} E_{\alpha, \alpha+\beta-1}(x) +  \frac{1-\beta}{\alpha} E_{\alpha, \alpha+\beta}(x)
\end{eqnarray}
all calculations may be reduced to Mittag-Leffler functions directly.

\begin{figure}
\begin{center}
\includegraphics[width=\textwidth]{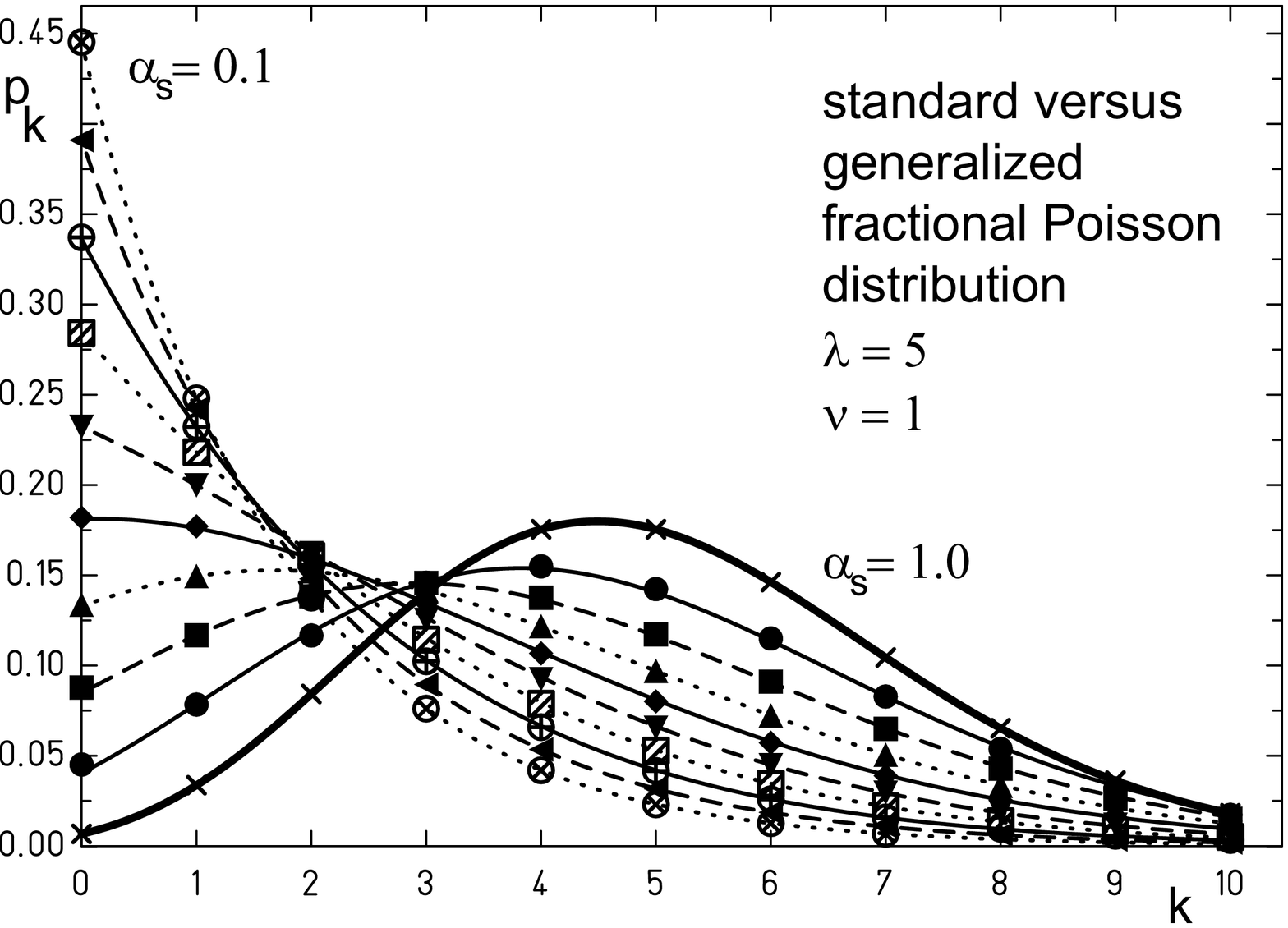}    
\caption{Fit results of the standard  fractional Poisson distribution $p_k (\alpha_s,\lambda)$ (\ref{poistart}) with $\nu=1$ for $0.1 \leq\alpha_s \leq 1.0$ in 0.1 steps (dots) with generalized fractional Poisson distribution from  $p_k (\alpha,\beta,\lambda)$ (\ref{poipnM}) (lines) for $\lambda=5$. Fit parameter sets $\{\alpha.\beta\}$ are listed in table  \ref{poistab1}.
}
\label{fig3}
\end{center}
\end{figure}

\section{Discussion}
\setcounter{equation}{0}
\setcounter{theorem}{0}
In contrast to the standard fractional Poisson distribution, which is valid only in the interval $0 \leq \alpha \leq 1$, our derivation
is not restricted to a specific region of allowed $\alpha$ values. Nevertheless, in the following we present results in the practically
interesting region $0 < \alpha  < 2$.
In figures \ref{fig1} and \ref{fig2} we have plotted the generalized Poisson distribution (\ref{poipnM}) for varying $\alpha$ and $\beta$ respectively.

To simplify an interpretation, we first calculate the mean and variance for two limiting cases:

First, for small $\lambda \ll 1$ we obtain
\begin{eqnarray}
\label{poiasymue}
\lim_{\lambda \rightarrow 0}\mu &=&
\frac{\Gamma(\beta)}{\Gamma(\alpha+\beta)} \lambda \\
\lim_{\lambda \rightarrow 0}\sigma &=&
\frac{\Gamma(\beta)}{\Gamma(\alpha+\beta)} \lambda 
\end{eqnarray}
which shows two interesting features:

The feature $\mu=\sigma$, which we know from the classical Poisson distribution, holds for
small $\lambda$ for the generalized case too.
For $\beta=1$ the mean corresponds to the mean value, which is obtained with standard fractional Poisson distribution (\ref{poistart}), only the shape differs (since variances are different).

However, this case ($\lambda\ll1$) is a very special case, in general, we deal with $\lambda\gg1$. 

For the case $\lambda \gg 1$, using the asymptotic expression for the generalized Mittag-Leffler func\-tion\cite{hau09}
\begin{eqnarray}
\label{poiasyminf}
\lim_{\lambda \rightarrow \infty}{E_{\alpha, \beta}(\lambda)} &=&
\frac{1}{\alpha}\lambda^{(1-\beta)/\alpha} e^{\lambda^{1/\alpha}} + o(\frac{1}{\lambda})
\qquad 0 <\alpha < 2
\end{eqnarray}
we obtain
\begin{eqnarray}
\label{poiasymue}
\lim_{\lambda \rightarrow \infty}\mu &=&
\frac{1-\beta + \lambda^{1/\alpha}}{\alpha} \\
\lim_{\lambda \rightarrow \infty}\sigma &=&
\frac{\lambda^{1/\alpha}}{\alpha^2}
\end{eqnarray}
which allows for a direct interpretation of parameter $\beta$: For a given $\alpha$, the expectation value $\mu$ is shifted by the amount $(1-\beta)/\alpha$, while the variance is independent of $\beta$. In figure \ref{fig2}, this behavior is easily observed.  

Consequently, the parameter $\beta$ could be directly interpreted within predator-prey systems: Since the probability distribution is shifted by a $\beta$ depended constant amount from the expected mean $\approx \lambda^{1/\alpha}/\alpha$,   while other characteristics (e.g. shape) of the distribution remain unchanged, $\beta$ may be considered as a constant time independent population decrease and increase rate respectively (comparable to a sink/source in electro dynamics).

With this understanding, in a next step we may then investigate the corresponding fractional process, considering a fractional pendant of the coupled sets of Lotka-Volterra equations\cite{mur08}, which are currently used to model such scenarios. 

Finally, the parameter set $\{\alpha,\beta\}$ in (\ref{poipnM}) may be adjusted appropriately to obtain similar  mean $\mu$ and variance $\sigma$ values of both, the standard fractional Poisson distribution (\ref{poistart}) and the generalized fractional Poisson distribution , which leads to similar shapes of the distribution function. 
This may be realized by calculating the parameter pair $\{\alpha,\beta\}$ by equating $\mu$ and $\sigma$ from 
(\ref{poimoments1}), (\ref{poimoments2}) with (\ref{poimoments12x}), (\ref{poimoments12y}), by the requirement, that both probability distributions should be same values for two given k's or by a fit procedure, e.g. near the maximum of the distribution.

In figure \ref{fig3} we compare the results of a fitting procedure of the fractional Poisson distribution (\ref{poistart}) with  (\ref{poipnM}). In table  \ref{poistab1} the corresponding fitted parameters are listed. The agreement is very good. In practical applications, it would be difficult, to distinguish both distribution families.

Therefore we may conclude, that the proposed generalized Poisson distribution covers the complete parameter range of the standard fractional Poisson distribution as a subset. 
It should at least be mentioned, that  the numerical behavior of the proposed generalized fractional Poisson distribution is much easier to handle than the standard distribution.

\begin{table}
\caption{
\label{poistab1}
Parameter sets for a fit of the standard fractional Poisson distribution $p_k (\alpha_s,\lambda)$  (\ref{poistart}) with $\nu=1$ for different $\alpha_s$ and the generalized fractional Poisson distribution  $p_k (\alpha,\beta,\lambda)$  (\ref{poipnM})  for $\lambda=5$.}
\begin{tabular}{l|ll||l|ll}
\hline\noalign{\smallskip}
$\alpha_s$& $\alpha$ & $\beta$&$\alpha_s$& $\alpha$ & $\beta$ \\ 
\noalign{\smallskip}\hline\noalign{\smallskip}
1.0& 1.000 & 1.00 &0.5&0.755&10.40\\ 
0.9& 0.911 & 2.86 &0.4&0.718&13.73\\ 
0.8& 0.855 & 4.58 &0.3&0.670&19.40\\ 
0.7& 0.815 & 6.33 &0.2&0.606&31.61\\ 
0.6& 0.781 & 8.26 &0.1&0.507&76.01\\ 
\hline\noalign{\smallskip}
\end{tabular}
\end{table} 

\section*{Conclusion}
We have presented a general method to derive the characteristic properties of a given fractional distribution in terms of normal ordered Bell polynomials with the Euler-operator as an argument.    

We have demonstrated for the special case of the generalized Mittag-Leffler function, that the investigation of corresponding  fractional distributions helps to obtain a clearer understanding of first the meaning of additional parameters in generalizations of Mittag-Leffler type functions  and second of possible areas of their application. 

In addition we have shown, that the standard fractional Poisson distribution may be
reproduced as a specific subset of the generalized fractional Poisson distribution. 

\section*{Acknowledgements}
We thank A. Friedrich for valuable discussions.




 \bigskip 
 
 \smallskip
 
 \it
 
 \noindent
$^1$ GigaHedron \\
Berliner Ring 80 \\
D-63303 Dreieich \\
GERMANY  \\[4pt]
e-mail: herrmann@gigahedron.com
\hfill Received: est. March 10, 2015 \\


\begin{thebibliography}{99}
 \normalsize

\bibitem{bel34}
E.~T. Bell,   (1934)
Exponential polynomials
\emph{Ann. Math.}  \textbf{35}(2) 258--277.

\bibitem{dob77}
G. Dobinski, (1877)
Summirung der Reihe $\sum n^m/n!$ f\"ur $m=1, 2, 3, 4, 5, ...$
\emph {Grunert Archiv (Arch. Math. Phys.)} \textbf{61} 333--336. 

\bibitem{hau09} 
H.~J. Haubold,  A.~M.  Mathai and  R.~K. Saxena,  (2009)
Mittag-Leffler functions and their applications
\emph{
arXiv:0909.0230, Journal of Applied Mathemathics, Hindawi},  (2011) 
298628.

\bibitem{kil13} 
A.~A. Kilbas, A.~A.  Koroleva and   S.~S. Rogosin,  (2013)
Multi-parameter Mittag-Leffler functions and their extension 
\emph{Fract. Calc. Appl. Anal.} \textbf{16} 378--404. 
DOI:10.2478/s13540-013-0024-9

\bibitem{kir99} 
V. Kiryakova, (1999)
Multi-indexed Mittag-Leffler functions, related Gelfond-Leontiev operators and Laplace type trnasforms
\emph{Fract. Calc. Appl. Anal.} \textbf{2} 445--462. 

\bibitem{las03} 
N. Laskin,  (2003) 
Fractional Poisson process
\emph{Commun. Nonlin. Sci. Num. Sim} \textbf{8} 201--213.
 

\bibitem{las09}
N. Laskin, (2009)
Some applications of the fractional Poisson probability distribution
\emph {J. Math. Phys.} \textbf{50} 113513. 


\bibitem{mee99} 
M.~M. Meerschaert, D.~A. Benson  and B. B\"aumer,
Multidimensional advection and fractional dispersion
\emph{Phys. Rev. E}  \textbf{59} (1999)  5026.




\bibitem{met00}  
R. Metzler and J. Klafter, 
The random walk's guide to anomalous diffusion: a fractional dynamics approach.
\emph{Phys. Rep.}  \textbf{339} (2000) 1--77;
DOI:10.1016/S0370-1573(00)00070-3


\bibitem{ml} 
 M.~G.  Mittag-Leffler, (1903) 
Sur la nouvelle function $E_{\alpha} (x)$
\emph{Comptes Rendus Acad. Sci. Paris} \textbf{137} 554--558.

\bibitem{mur08}  
J.~D. Murray (2008)
Mathematical Biology I: An Introduction. 3. ed. 
Springer, Berlin

\bibitem{pod99}  
I. Podlubny, (1999)
\emph{Fractional Differential Equations}
Academic Press, Boston

\bibitem{rep00} 
O.~N. Repin  and A.~I. Saichev,    (2000) 
Fractional Poisson law
\emph{Radiophys. Quant. Electron.} \textbf{43} 738--741.

\bibitem{rom84} 
S. Roman, (1984)
\emph{ "The Exponential Polynomials" and "The Bell Polynomials." §4.1.3 and §4.1.8}
 in The Umbral Calculus. New York: Academic Press, pp. 63--67 and 82--87. 

\bibitem{sam93} 
S.~G. Samko,  A.~A. Kilbas   and  O.~I. Marichev,  
\emph{Fractional integrals and derivatives} 
Translated from the 1987 Russian original, Gordon and Breach, Yverdon(1993)

\bibitem{sti30}
J.  Stirling, (1730)
Methodus differentialis, sive tractatus de summation et interpolation serierum infinitarium. London. 
English translation by Holliday, J. The Differential Method: A Treatise of the Summation and Interpolation of Infinite Series. 1749.  

\bibitem{wic50} 
G.~C. Wick,   (1950)
The evaluation of the collision matrix
\emph{Phys. Rev.} \textbf{80} 268.

\bibitem{wi} 
A. Wiman,  (1905) 
\"Uber den Fundamentalsatz in der Theorie der Funktionen $E_a (x)$
\emph{Acta Math.} \textbf{29} 191--201.



\end{thebibliography}
\end{document}